\documentclass[final,1p,times]{elsarticle}
\usepackage{lmodern}
\usepackage{amssymb}
\usepackage{bm}
\usepackage{graphicx,graphics}
\usepackage{caption2}
\usepackage{psfrag}
\usepackage{amsmath,amsfonts}
\usepackage{booktabs}
\usepackage{amscd}
\usepackage{multirow}
\usepackage{lmodern}
\usepackage{float}
\usepackage{latexsym}
\usepackage{lscape}
\usepackage{color}

\newtheorem{theorem}{Theorem}[section]

\newcommand{\ds}{\displaystyle}

\newcommand{\om}{\Omega}
\newcommand{\p}{\partial}

\journal{}

\begin{document}

\bibliographystyle{plain}
\begin{frontmatter}

\title{Analysis of evolution equation with
variable-exponent memory modeling multiscale viscoelasticity
\footnote{Corresponding author: Xiangcheng Zheng}}

\author[1]{Yiqun Li}
\ead{YiqunLi24@outlook.com}

\author[2]{Xiangcheng Zheng}
\ead{xzheng@sdu.edu.cn}

\address[1]{School of Mathematics and Statistics, Wuhan University,
Wuhan 430072, China}

\address[2]{School of Mathematics, State Key Laboratory of Cryptography and Digital Economy Security, Shandong University, Jinan 250100, China}

\begin{abstract}
We investigate the well-posedness and solution regularity of an evolution equation with non-positive type variable-exponent memory, which describes multiscale viscoelasticity in materials with memory. The perturbation method is applied for model transformation, based on which the well-posedness is proved. Then the weighted solution regularity is derived, where the initial singularity is characterized by the initial value of variable exponent.
\end{abstract}


\end{frontmatter}

\section{Introduction}
This work considers the mathematical analysis of  the evolution equation with non-positive type variable-exponent memory, where the variable exponent  determined by the fractal dimension of the microstructures of viscoelastic materials describes the structure variations under the long-term cyclic loads, which in turn propagate to macro scales that eventually results in material failure \cite{Dan,FanHu,GarGiu,HouLow,HouWu,LiWan,LiaSty}
     \begin{equation}\label{VtFDEs}
     \begin{split}
         &~~~\p_t u(\bm x,t) -k(t)* \Delta u(\bm x,t) = f(\bm x,t), \quad (\bm x,t)\in \Omega\times (0,T];\\[0.05in]
         & \ds u(\bm x,0)=u_0(\bm x),\quad\bm x\in\Omega;\quad u(\bm x,t)=0,\quad (\bm x,t)\in \partial\Omega\times [0,T].
     \end{split}
     \end{equation}
    Here $\Omega\subset \mathbb R^d$ ($1\leq d \leq 3$) is a simply-connected bounded domain with the piecewise smooth boundary $\partial\Omega$ and convex corners, $T>0$, the source term $f(\bm x,t)$ and $u_0(\bm x)$ are given functions, and the convolution term is defined via the variable-exponent kernel $k(t)$ with $\alpha(t) \in (0,1)$
     \begin{equation}\label{k}
     \begin{split}
      k(t)*\Delta u(\bm x,t):=\int_0^tk(t-s)\Delta u(\bm x,s)ds,~~     k(t) : = \frac{t^{-\alpha(t)}}{\Gamma(1-\alpha(t))}.
     \end{split}
     \end{equation}

    The mathematical analysis for model \eqref{VtFDEs} with $\alpha(t) \equiv \hat \alpha$ for some $0<\hat \alpha <1$ has attracted extensive research activities, see e.g. \cite{MclTho,Da,YiGuo}, while, to our best knowledge, the variable-exponent  model \eqref{VtFDEs} remains untreated in the literature since the variable-exponent kernel $k(t)$ loses favorable properties  of its constant-exponent counterpart such as positive definiteness. Consequently, many established techniques for constant-exponent models do not apply.

    Recently, a perturbation method is developed in \cite{Zhe} to analyze the variable-exponent two-sided space-fractional problems, which splits the variable-exponent kernel \eqref{k} into its constant-order analogue and a low-order perturbation term. We adopt this approach to reformulate the model \eqref{ode1} to a more tractable form, based on which the well-posedness of the model \eqref{VtFDEs} and its regularity estimate are proved. Different from the resolvent estimate method in \cite{Zhe}, we employ special functions in solution expression for more subtle estimates. The developed results provide a theoretical support for numerical analysis and real applications of the proposed model. 

\section{Preliminaries}\label{sec2}
\subsection{Notations}
Let $L^p(\om)$ with $1 \le p \le \infty$ be the Banach space of $p$th power Lebesgue integrable functions on $\om$. For a positive integer $m$,
let  $ W^{m, p}(\Omega)$ be the Sobolev space of $L^p$ functions with $m$th weakly derivatives in $L^p(\om)$ (similarly defined with $\om$ replaced by an interval $\mathcal I$). Let  $H^m(\Omega) := W^{m,2}(\Omega)$ and $H^m_0(\Omega)$ be its subspace with the zero boundary condition up to order $m-1$. 
For a Banach space $\mathcal{X}$, let $W^{m, p}(0,T; \mathcal{X})$ be the space of functions in $W^{m, p}(0,T)$ with respect to $\|\cdot\|_{\mathcal {X}}$. All spaces are equipped with standard norms \cite{AdaFou}.

 Denote the eigenpairs of $-\Delta$ with Dirichlet boundary conditions by $\{\lambda_i,\phi_i\}_{i=1}^{\infty}$ in which $\{\phi_i\}_{i=1}^{\infty}$ form an orthonormal basis in $L^2(\om)$ and the eigenvalues $\{\lambda_i\}_{i=1}^\infty$ form a positive and  non-decreasing sequence \cite{Eva}.
  We introduce the Sobolev space $\check{H}^s(\Omega)$ for $s\geq 0$ by
$$ \check{H}^{s}(\Omega) := \big \{ q \in L^2(\Omega): \| q \|_{\check{H}^s(\Omega)}^2 : = \sum_{i=1}^{\infty} \lambda_i^{s} (q,\phi_i)^2 < \infty \big \},$$
which is a subspace of $H^s(\Omega)$ satisfying $\check{H}^0(\Omega) = L^2(\Omega)$ and $\check{H}^2(\Omega) = H^2(\Omega) \cap H^1_0(\Omega)$ \cite{Tho}.

Through this paper, we denote $\alpha_0:= \alpha(0)$ and assume that $\alpha(t)$ has bounded second-order derivative on $[0, T]$.  In addition, we use $Q$ to denote a generic positive constant that may assume different values at different occurrences. For simplicity, we may drop the domain $\om$ and the time interval $(0,T)$ in the Sobolev spaces and norms, e.g. we write $\| \cdot\|_{L^2(L^2)}$ instead of $\| \cdot\|_{L^2(0,T;L^2(\Omega))}$, when no confusion occurs.

\subsection{Spectral decomposition}
For $t \in [0, T]$, we expand $u$ and $f$ with respect to  $\{\phi_i\}_{i=1}^{\infty}$ as follows \cite{SakYam}
$$u=\sum_{i=1}^\infty u_i(t)\phi_i(\bm x),~u_i(t) := \big (u(\cdot,t),\phi_i \big),~~f=\sum_{i=1}^\infty f_i(t)\phi_i(\bm x),~f_i(t) := \big (f(\cdot,t),\phi_i \big).$$
We invoke these expressions into (\ref{VtFDEs}) to find that $\{u_i\}_{i=1}^{\infty}$ satisfy the following integro-differential equations
\begin{align}\label{ode1}
u_i'+\lambda_i k*u_i=f_i,\text{ for }t\in (0,T];~~u_i(0)=u_{0,i}:=(u_0,\phi_i), \quad i=1,2, \cdots
\end{align}

 For $\mu>0$, we first define
$\beta_\mu (t):=\frac{t^{\mu-1}}{\Gamma(\mu)}$.  We follow the idea of the perturbation method proposed in \cite{Zhe} to split the kernel $k(t)$ as
\begin{align}
k(t)&=\frac{t^{-\alpha(t)}}{\Gamma(1-\alpha(t))}=\frac{t^{-\alpha_0}}{\Gamma(1-\alpha_0)}+\int_0^t\p_z\frac{t^{-\alpha(z)}}{\Gamma(1-\alpha(z))}dz=:\beta_{1-\alpha_0}+\tilde g(t).\label{hh3}
\end{align}
It is shown in \cite{Zhe} that
\begin{align}\label{dbndg}
|\tilde g|\leq Q\int_0^t t^{-\alpha_0}(1+|\ln t|)dz\leq Qt^{1-\alpha_0}(1+|\ln t|),~~|\tilde g'|\leq Qt^{-\alpha_0}(1+|\ln t|).
\end{align}
We thus employ (\ref{hh3}) to rewrite (\ref{ode1}) as
\begin{align}\label{ode2}
u_i'+\lambda_i \beta_{1-\alpha_0}*u_i=f_i-\lambda_i\tilde g*u_i,\text{ for }t\in (0,T];~~u_i(0)=u_{0,i}.
\end{align}

\section{Mathematical analysis}\label{sec3}

We prove the well-posedness and the regularity results  of (\ref{VtFDEs}).

\begin{theorem}\label{thm1:pde}
Suppose  $f \in H^1(L^2)$ and $u_{0} \in \check H^2$, the problem (\ref{VtFDEs})  has a unique solution $u \in H^1(L^2)$ and
\begin{equation}\label{pdestab1}\begin{array}{l}
\ds \|u\|_{H^1(L^2)}\leq Q\big(\|u_{0}\|_{\check H^2}+ \|f\|_{H^1(L^2)}\big).
\end{array}
\end{equation}
\end{theorem}
\textbf{Proof}
The proof of the theorem  will be carried out in two steps.

\textit{Step A: Well-posedness of problem (\ref{ode2}).}

We first consider the case $u_{0,i}=0$ in (\ref{ode2}) and let $\mathcal X := \{ q \in H^1(0,T): q(0)=0 \}$ equipped with the equivalent norm $\|q\|_{\mathcal X,\sigma} := \|e^{-\sigma t} q'\|_{L^2(0,T)}$ for some $\sigma\geq 0$. For each $v\in \mathcal X$, let $w:= \mathcal M v$ be the solution of
\begin{align}\label{ode4}
w'+\lambda_i \beta_{1-\alpha_0}*w=f_i-\lambda_i\tilde g*v,\text{ for }t\in (0,T];~~w(0)=0.
\end{align}
Then  the solution $w$ to (\ref{ode4}) could be represented as \cite{Mcl}
\begin{align}\label{OdeSol}
w=E_{2-\alpha_0,1}(-\lambda_it^{2-\alpha_0})*(f_i-\lambda_i\tilde g*v).
\end{align}
To bound $\|w\|_{\mathcal X,\sigma}$, we directly differentiate the above equation to get
\begin{align}\label{w'}
w'=E_{2-\alpha_0,1}(-\lambda_it^{2-\alpha_0})f_i(0)+E_{2-\alpha_0,1}(-\lambda_it^{2-\alpha_0})*(f_i'-\lambda_i\tilde g'*v).
\end{align}
By using the Mittag-Leffler function \cite{DieFor,Gor,Pod} and the fact that $$[tE_{2-\alpha_0,2}(-\lambda_it^{2-\alpha_0})\big]_t = E_{2-\alpha_0,1}(-\lambda_it^{2-\alpha_0}),$$ we further apply integration by parts to reformulate the last term on the right-hand side of (\ref{w'}) to get
\begin{equation}\label{ode5}\begin{array}{rl}
\hspace{-0.1in} \ds E_{2-\alpha_0,1}(-\lambda_it^{2-\alpha_0})*(\lambda_i\tilde g'*v)&\ds =tE_{2-\alpha_0,2}(-\lambda_it^{2-\alpha_0})*(\lambda_i\tilde g'*v')\\[0.1in]
\hspace{-0.1in} \ds &\ds =[t^{\alpha_0-1}(\lambda_i t^{2-\alpha_0})E_{2-\alpha_0,2}(-\lambda_it^{2-\alpha_0})]*(\tilde g'*v').
\end{array}
\end{equation}
By the asymptotic property of the Mittag-Leffler function \cite{Gor,Jinbook},  we have
\begin{align}\label{MLbound}
\big|(\lambda_i t^{2-\alpha_0})E_{2-\alpha_0,2}(-\lambda_it^{2-\alpha_0})\big|\leq Q,
\end{align}
where $Q$ is independent from $t$ and $\lambda_i$. Thus we apply Young's convolution inequality and the estimate
\begin{equation}\label{AdjEq:e4}
\int_0^Te^{-\sigma t }t^{\alpha_0-1}dt\leq \sigma^{-\alpha_0}\int_0^{\infty}e^{-z}z^{\alpha_0-1}dz\leq Q\sigma^{-\alpha_0},
\end{equation}
and (\ref{MLbound}) to bound (\ref{ode5}) as
\begin{equation*}\begin{array}{rl}
&\big\|e^{-\sigma t}[E_{2-\alpha_0,1}(-\lambda_it^{2-\alpha_0})*(\lambda_i\tilde g'*v)]\big\|_{L^2(0,T)}\\[0.1in]
& \qquad\ds\leq Q\big\|e^{-\sigma t}(t^{\alpha_0-1}*|\tilde g'*v'|)\big\|_{L^2(0,T)}\\[0.1in]
&\qquad\ds\leq Q\big\|(e^{-\sigma t}t^{\alpha_0-1})*|e^{-\sigma t}\tilde g'|*|e^{-\sigma t}v'|\big\|_{L^2(0,T)}\\[0.1in]
&\qquad \ds\leq Q\|e^{-\sigma t}t^{\alpha_0-1}\|_{L^1(0,T)}\|e^{-\sigma t}\tilde g'\|_{L^1(0,T)}\|v\|_{\mathcal X,\sigma}\leq Q\sigma^{-\alpha_0}\|v\|_{\mathcal X,\sigma}.
\end{array}
\end{equation*}
We invoke this estimate, the boundedness of $E_{2-\alpha_0,1}(-\lambda_it^{2-\alpha_0})$ \cite{Gor,Jinbook} and $|f_i(0)|\leq Q\|f_i\|_{H^1(0,T)}$ in (\ref{w'}) to get
\begin{align}\label{ode7}
\|w\|_{\mathcal X,\sigma}\leq Q\|f_i\|_{H^1(0,T)}+Q\sigma^{-\alpha_0}\|v\|_{\mathcal X,\sigma}.
\end{align}
Thus the mapping $\mathcal M:\mathcal X\rightarrow \mathcal X$ is well defined. To show its contractivity, let $w_i=\mathcal M v_i$ for $i =1$, $2$,  $e_w : = w_1- w_2$ and $e_v : = v_1- v_2$ satisfy
\begin{align}\label{ode6}
e_w'+\lambda_i \beta_{1-\alpha_0}*e_w=-\lambda_i\tilde g*e_v.
\end{align}
 Then an application of the estimate (\ref{ode7}) to (\ref{ode6})  gives
 \begin{align}\label{ode8}
\|e_w\|_{\mathcal X,\sigma}\leq Q\sigma^{-\alpha_0}\|e_v\|_{\mathcal X,\sigma}.
\end{align}
 We select a sufficiently large $\sigma$ such that the mapping $\mathcal M$ is a contraction. By the Banach fixed point theorem, $\mathcal M$ has a unique fixed point $v=\mathcal Mv$ with the estimate $\|v\|_{\mathcal X,\sigma}\leq Q\|f_i\|_{H^1(0,T)}$.

For the case that $u_{0,i}\neq 0$, we could directly verify that
$$u_i=E_{2-\alpha_0,1}(-\lambda_it^{2-\alpha_0})u_{0,i}+v,$$
where $v$ is the fixed point of $\mathcal M$, solves
(\ref{ode2}) and thus (\ref{ode1}). Differentiate this equation to get
$$u_i'=-\lambda_i t^{1-\alpha_0}E_{2-\alpha_0,2-\alpha_0}(-\lambda_it^{2-\alpha_0})u_{0,i}+v',$$
which, together with $\|v\|_{\mathcal X,\sigma}\leq Q\|f_i\|_{H^1(0,T)}$, leads to
\begin{equation}\label{odestab}\begin{array}{l}
\ds \|u_i\|_{H^1(0,T)}\leq Q\big(\lambda_i|u_{0,i}|+\|f_i\|_{H^1(0,T)}\big).
\end{array}
\end{equation}
The uniqueness of (\ref{ode1}) follows from that of (\ref{ode2}) with $u_{0,i}=0$. Consequently, we conclude that  (\ref{ode2}) and thus (\ref{ode1}) admit  a unique solution in $H^1(0,T)$ with the estimate (\ref{odestab}).

\textit{Step B: Well-posedness of the problem  (\ref{VtFDEs}).}

We invoke the  estimate   (\ref{odestab}) to bound
 $   \bar u:=\int_0^t\big[\sum_{i=1}^\infty u_i'(s)\phi_i(\bm x)\big]ds + u_0
$
by using
\begin{equation}\begin{array}{rl}\label{AdjEq:e8}
\ds \|\bar u\|_{H^1(L^2)}^2 \leq Q \|\p_t \bar u\|^2_{L^2(L^2)} = \sum_{i=1}^\infty \|u_i'\|_{L^2(0,T)}^2
&\hspace{-0.1in}\ds \leq  Q\sum_{i=1}^\infty \big(\lambda_i^2|u_{0,i}|^2 + \|f_i\|_{H^1(0,T)}^2\big)\\[0.15in]
&\hspace{-0.1in}\ds  =Q\big(\|u_0\|_{\check H^2}^2 + \|f\|^2_{H^1(L^2)}\big).
\end{array}
\end{equation}
As$u_i(t) = \int_0^t u_i'(s)ds + u_{0,i}$ satisfies the differential equation in (\ref{ode1}) for $i \ge 1$, $\bar u \in H^1(L^2)$ is a solution to problem (\ref{VtFDEs}).
The uniqueness of the solution to the problem (\ref{VtFDEs})   follows from that of the differential equation in (\ref{ode1}). We thus complete the proof of the theorem.

\begin{theorem}\label{thm:pde}
Suppose  $f \in H^1(\check H^2)$ and $u_{0} \in \check H^4$, the following regularity estimate holds for the  problem \eqref{VtFDEs}
\begin{equation}\label{pde:reg}\begin{array}{l}
\ds  \big\|t^{\alpha_0/2} \p_t^2 u \big\|_{L^2(L^2)}  +  \| u  \|_{H^1(\check H^2)} \le Q\big(\|u_{0}\|_{\check H^4}+ \|f\|_{H^1(\check H^2)}\big).
\end{array}
\end{equation}
\end{theorem}
\textbf{Proof}
We  differentiate the equation in (\ref{ode2}) to obtain
 \begin{align}\label{odeutt}
u_i^{\prime \prime}+\lambda_i \beta_{1-\alpha_0}*u_i^{\prime}=f_i^{\prime}-\lambda_i\tilde g^{\prime}*u_i-\lambda_i \beta_{1-\alpha_0}u_{0,i},
\end{align}
which, combined with the first estimate in (\ref{odestab}), (\ref{dbndg}) and Young’s convolution
inequality, gives
\begin{equation}\label{utt}\begin{array}{rl}
\ds \big\|t^{\alpha_0/2} u_i^{\prime \prime}\big\|_{L^2(0, T)} &\hspace{-0.1in}\ds  \le Q \big(\lambda_i \|\beta_{1-\alpha_0}\|_{L^1(0, T)}\|u_i^{\prime}\|_{L^2(0, T)} + \|f_i^{\prime}\|_{L^2(0, T)} \\[0.1in]
&\hspace{-0.1in}\ds  \qquad + \lambda_i \|\tilde g^{\prime}\|_{L^1(0, T)}\|u_i\|_{L^2(0, T)} + \lambda_i|u_{0,i}| \big)\\[0.1in]
&\hspace{-0.1in}\ds \le Q\big(\lambda_i\|u_i\|_{H^1(0, T)} + \|f_i^{\prime}\|_{L^2(0, T)} + \lambda_i|u_{0,i}|  \big)\\[0.1in]
&\hspace{-0.1in}\ds   \le Q\big(\lambda_i^2|u_{0,i}|+ \lambda_i \|f_i\|_{H^1(0,T)}\big).
\end{array}
\end{equation}
We further combine   (\ref{odestab}) and (\ref{utt}) and follow the procedures in (\ref{AdjEq:e8})  to obtain
\begin{equation*}\begin{array}{l}
\ds \big \|t^{{\alpha_0}/{2}}\p_t^2 u \big \|^2_{L^2(L^2)} +  \| u  \|^2_{H^1(\check H^2)}  = \sum_{i=1}^\infty \big( \|t^{{\alpha_0}/{2}} u_i^{\prime \prime}\|_{L^2(0,T)}^2 + \lambda_i^2\| u_i\|_{H^1(0,T)}^2 \big)\\[0.1in]
\qquad \qquad \ds \leq Q\sum_{i=1}^\infty \big( \lambda_i^2\|f_i\|^2_{H^1(0,T)} + \lambda_i^4 |u_{0,i}|^2 \big)\le Q \big ( \|f\|^2_{H^1(\check H^2)} + \|u_0\|^2_{\check H^4} \big),
\end{array}\end{equation*}
which gives (\ref{pde:reg}) and we thus complete the proof.

\section{Concluding remarks}

This work establishes the well-posedness and regularity results of an evolution equation with non-positive type variable-exponent memory.
A perturbation method is employed to reformulate the model, based on which the analysis is performed. Based on these theoretical results, we may consider error estimates of the numerical approximation to model \eqref{VtFDEs}. However, the coupling of the variable-exponent kernel and the Laplacian may cause difficulties in  stability and error estimates that have not been encountered before. We will investigate this interesting topic in the near future.

\section*{Acknowledgments}
This work was partially supported by the National Key R\&D Program of China (No. 2023YFA1008903), the National Natural Science Foundation of China (No. 12301555), the Taishan Scholars Program of Shandong Province (No. tsqn202 306083),  the China Postdoctoral Science Foundation (No. 2024M762459) and the Natural Science Foundation of Hubei Province (No. JCZRQN202500278).

\end{document}